\documentclass[a4paper]{article}
\pdfoutput=1

\usepackage{graphicx,wrapfig}
\usepackage{lipsum}
\usepackage{amsmath,amssymb,mathrsfs}
\usepackage{hyperref}
\usepackage{url}
\usepackage{mathtools}
\usepackage{changepage}
\usepackage{multirow}
\usepackage{wrapfig}
\usepackage{subfig}
\usepackage{color}
\usepackage{epsfig,enumerate,amsmath,amsfonts,amssymb,amsthm,mathrsfs,ifpdf}
\usepackage{indentfirst,relsize}
\usepackage{setspace,graphicx}
\usepackage{latexsym}
\usepackage[all]{xy}
\usepackage[usenames,dvipsnames]{pstricks}
\usepackage{pst-grad} % For gradients
\usepackage{pst-plot} % For axes
\usepackage[margin = 2.50cm]{geometry}

\usepackage{algorithm}

\usepackage{algpseudocode}

\newcommand{\remove}[1]{}

\newtheorem{theo}{Theorem}%[section]
\newtheorem{lem}[theo]{Lemma}

\newtheorem{cl}[theo]{Claim}
\newcounter{Ca}[theo]
\newtheorem{ca}[Ca]{Case}

\title{Strong transitivity of a graph}
\author{Subhabrata Paul \and Kamal Santra}
\author{Subhabrata Paul\footnote{Department of Mathematics, IIT Patna, India, email:subhabrata@iitp.ac.in} \and Kamal Santra\footnote{Department of Mathematics, IIT Patna, India, email:kamal\_1821ma04@iitp.ac.in} }

\date{}

\begin{document}

\maketitle
\begin{abstract}
	A vertex partition $\pi = \{V_1, V_2, \ldots, V_k\}$ of $G$ is called a \emph{transitive partition} of size $k$ if $V_i$ dominates $V_j$ for all $1\leq i<j\leq k$. For two disjoint subsets $A$ and $B$ of $V$, we say  $A$ \emph{strongly dominates} $B$ if for every vertex $y\in B$, there exists a vertex $x\in A$, such that $xy\in E$ and $deg_G(x)\geq deg_G(y)$. A vertex partition $\pi = \{V_1, V_2, \ldots, V_k\}$ of $G$ is called a \emph{strong transitive partition} of size $k$ if $V_i$ strongly dominates $V_j$ for all $1\leq i<j\leq k$. The \textsc{Maximum Strong Transitivity Problem} is to find a strong transitive partition of a given graph with the maximum number of parts. In this article, we initiate the study of this variation of transitive partition from algorithmic point of view. We show that the decision version of this problem is NP-complete for chordal graphs. On the positive side, we prove that this problem can be solved in linear time for trees and split graphs.

\end{abstract}

{\bf Keywords.}
Strong transitivity, NP-completeness, Linear-time algorithm, Trees, Split graphs, Chordal graphs.
%%%%%%%%%%%%%%%%%%%%%%%%%%%%%%%%%%%%%%%%%%%%%%%%%%%%%%%%%%%%%%%%%%%%%%%%%%%%%%%%%%%%%%%%%%%%%%%%%%%%%%
\section{Introduction}
Partitioning a graph is one of the fundamental problems in graph theory. In the partitioning problem, the objective is to partition the vertex set (or edge set) into some parts with desired properties, such as independence, minimal edges across partite sets, etc. A \emph{dominating set} of $G=(V, E)$ is a subset of vertices $D$ such that every vertex $x\in V\setminus D$ has a neighbour $y\in D$, that is, $x$ is dominated by some vertex $y$ of $D$. For two disjoint subsets $A$ and $B$ of $V$, we say $A$ \emph{dominates} $B$ if every vertex of $B$ is adjacent to at least one vertex of $A$. Many variants of partitioning problem have been studied in literature based on some domination relationship among the partite sets. For example \emph{domatic partition}\cite{cockayne1977towards,zelinka1980domatically} (each partite set is a dominating set), \emph{Grundy partition}\cite{hedetniemi1982linear,zaker2005grundy,zaker2006results} (each partite set is independent and dominates every other partite sets after itself), \emph{transitive partition} \cite{hedetniemi2018transitivity,haynes2019transitivity,paul2023transitivity,santra2023transitivity} (a generalization of Grundy partition where partite sets need not be independent), \emph{upper domatic partition} \cite{haynes2020upper,samuel2020new} (a generalization of transitive partition where for any two partite sets $X$ and $Y$ either $X$ dominates $Y$ or $Y$ dominates $X$ or both).

In 1996, Sampathkumar and Pushpa Latha introduced the notion of \emph{strong domination} \cite{sampathkumar1996strong}. A \emph{strong dominating set} of $G=(V, E)$ is a subset of vertices $D$ such that for every vertex $x\in V\setminus D$, $x$ is dominated by some vertex $y\in D$ and $deg_G(y)\geq deg_G(x)$. Recently, based on this strong domination, a variation of domatic partition, namely \emph{strong domatic partition}, has been studied in \cite{ghanbari2023strong}. In the \emph{strong domatic partition}, the vertex set is partitioned into $k$ parts, say $\pi =\{V_1, V_2, \ldots, V_k\}$, such that each $V_i$ is a strong dominating set of $G$. In this article, we introduce a variation of transitive partition based on strong domination, namely \emph{strong transitive partition}. A vertex partition $\pi = \{V_1, V_2, \ldots, V_k\}$ of $G$ is called a \emph{strong transitive partition} of size $k$ if $V_i$ strongly dominates $V_j$ for all $1\leq i<j\leq k$. The maximum order of such a strong transitive partition is called \emph{strong transitivity} of $G$ and is denoted by $Tr_{st}(G)$. The \textsc{Maximum Strong Transitivity Problem} and its corresponding decision version are defined as follows:

\noindent\textsc{\underline{Maximum Strong Transitivity Problem(MSTP)}}

\noindent\emph{Instance:} A graph $G=(V,E)$

\noindent\emph{Solution:} A strong transitive partition of $G$

\noindent\emph{Measure:} Order of the strong transitive partition of $G$

\noindent\textsc{\underline{Maximum Strong Transitivity Decision Problem(MSTDP)}}

\noindent\emph{Instance:} A graph $G=(V,E)$, integer $k$

\noindent\emph{Question:} Does $G$ have a Strong transitive partition of order at least $k$?

Note that every strong transitive partition is also a transitive partition. Therefore, for any graph $G$, $1\leq Tr_{st}(G)\leq Tr(G)\leq \Delta(G)+1$, where $\Delta(G)$ is the maximum degree of $G$. From the definition of a strong transitive partition, it is clear that for the regular graph, transitivity is same as strong transitivity; as a consequence, for the graph class $K_n$ and cycle $C_n$, they are the same. However, the transitive partition of a graph is not always a strong transitive partition, even for the same value of both parameters. For a path $P_3$, with vertex set $\{a, b, c\}$, tacking $\pi=\{V_1={a, c}, V_2=\{b\}\}$, then it is a transitive partition but not a strong transitive partition as the $deg(b)>deg(a)$ and $deg(b)>deg(a)$. But considering $\pi'=\{V_1={b, c}, V_2=\{a\}\}$, then it is both a strong transitive and transitive partition with size $2$. It can be easily verified that for the graph class $P_n, n\geq 6$, transitivity and strong transitivity are the same and equal to $3$. So, we see that there are graph classes where both parameters have the same value, but generally, their difference can be arbitrarily large. If $G$ is a complete bipartite graph of the form $K_{m, m-1}$, then $Tr_{st}(K_{m, m-1})=2$. Let $V(G)=X\cup Y$. Also, let $x\in X$ and consider a vertex partition $\pi=\{V_1, V_2\}$, where $V_1=(X\setminus\{x\})\cup Y$, $V_2=\{x\}$. Since $m\geq2$, there exits $y\in Y$ and $m=deg(y)\geq deg(x)=m-1$. So, $\pi$ is a strong transitive partition of $G$. Therefore, $Tr_{st}(K_{m, m-1})\geq 2$. To prove $Tr_{st}(K_{m, m-1})=2$, we now show $Tr_{st}(K_{m, m-1})<3$ by contradiction. Assume $Tr_{st}(K_{m, m-1})\geq 3$ and $\pi=\{V_1, V_2, \ldots, V_k\}$ be a strong transitive partition of $G$ of size $k$. Let $y\in V_i$, $3\leq i\leq k$ and $y\in Y$. Since $\pi$ is a strong transitive partition, $V_1$ strongly dominates $V_i$. So, for $y\in V_i$, there must exists a vertex $x\in V_1$, such that $xy\in E(G)$ and $deg_G(x)\geq deg_G(y)$. But for every vertex $x\in X$, $deg_G(x)=m-1< deg_G(y)=m$. So, $V_i$ cannot contain vertices from $Y$. Therefore, $V_i$ contains only vertices from $X$. Let $x\in V_i$ and $x\in X$, $3\leq i\leq k$. As $\pi$ is a strong transitive partition, $V_2$ strongly dominates $V_i$. So, for $x\in V_i$, there must exists a vertex $y\in V_2$, such that $xy\in E(G)$ and $deg_G(y)\geq deg_G(x)$. Now, $V_1$ strongly dominates $V_2$ and $y\in V_2$. To strongly dominates $y$, we have a vertex $x'\in V_1\cup X$, such that $x'y\in E(G)$ and $deg_G(x')\geq deg_G(y)$. Again, this is not possible as the $deg_G(x')=m-1< deg_G(y)=m$. Therefore, $V_i$, $3\leq i\leq k$ cannot contain vertices from $X$ also. So, $k\leq 3$. Therefore, if $G$ is a $K_{m, m-1}$, $m\geq 2$, then $Tr_{st}(K_{m, m-1})=2$. We know that $Tr(K_{m, m-1})=\min\{m+1, m\}=m$, so we have $Tr_{st}(K_{m, m-1})=2$. So, the difference $Tr(G)-Tr_{st}(G)=m-2$, which is arbitrarily large for $m$. For the transitivity, if $H$ is a subgrah of graph $G$, then $Tr(H)\leq Tr(G)$\cite{hedetniemi2018transitivity}. But for the strong transitivity, this is not true. Considering $G=K_{3, 2}$ and $H=C_4$. Clearly, $H$ is a subgraph of $G$, and $Tr_{st}(G=K_{3, 2})=2$, $Tr_{st}(C_4)=3$. So, in this example, $Tr_{st}(H)>Tr_{st}(G)$. Moreover, the behaviour of the strong transitivity is the same as the transitivity when the graph is disconnected. It is also known that every connected graph $G$ with $Tr(G)=k\geq 3$ has a transitive partition $\pi =\{V_1,V_2, \ldots, V_k\}$ such that $|V_k|$ = $|V_{k-1}| = 1$ and $|V_{k-i}| \leq 2^{i-1}$ for $2\leq i\leq k-2$ \cite{haynes2019transitivity}. This implies that the maximum transitivity problem is fixed-parameter tractable \cite{haynes2019transitivity}. Since a strong transitive partition of a graph is also a transitive partition, MSTP is also fixed-parameter tractable.

In this paper, we study the computational complexity of this problem. The main contributions are summarized below:
\begin{enumerate}
	\item [1.] The \textsc{MSTDP} is NP-complete for chordal graphs.
	
	\item[2.] The \textsc{MSTP} can be solved in linear time for trees and split graphs.
	
\end{enumerate}
The rest of the paper is organized as follows. Section 2 shows that the \textsc{MSTDP} is NP-complete in chordal graphs. Section 3 describes linear-time algorithms for trees and split graphs. Finally, Section 4 concludes the article.

\section{NP-complete for chordal graphs of strong transitivity}
This section shows that \textsc{Maximum Strong Transitivity Decision Problem} is NP-complete for chordal graphs. A graph is called \emph{chordal} if there is no induced cycle of length more than $3$. Clearly, \textsc{MSTDP} is in NP. We prove the NP-completeness of this problem by showing a polynomial-time reduction from \textsc{Proper $3$-Coloring Decision Problem} in graphs, which is known to be NP-complete \cite{garey1990guide}. A proper $3$-colring of a graph $G=(V,E)$ is a function $g$, from $V$ to $\{1,2,3\}$, such that for any edge $uv \in E$, $g(u)\not= g(v)$. The  \textsc{Proper $3$-Coloring Decision Problem} is defined as follows:

\noindent\textsc{\underline{Proper $3$-Coloring Decision Problem (P$3$CDP)}}

\noindent\emph{Instance:} A graph $G=(V,E)$

\noindent\emph{Question:} Does there exist a proper $3$-coloring of $G$?

\begin{figure}[htbp!]
	\centering
	\includegraphics[scale=0.45]{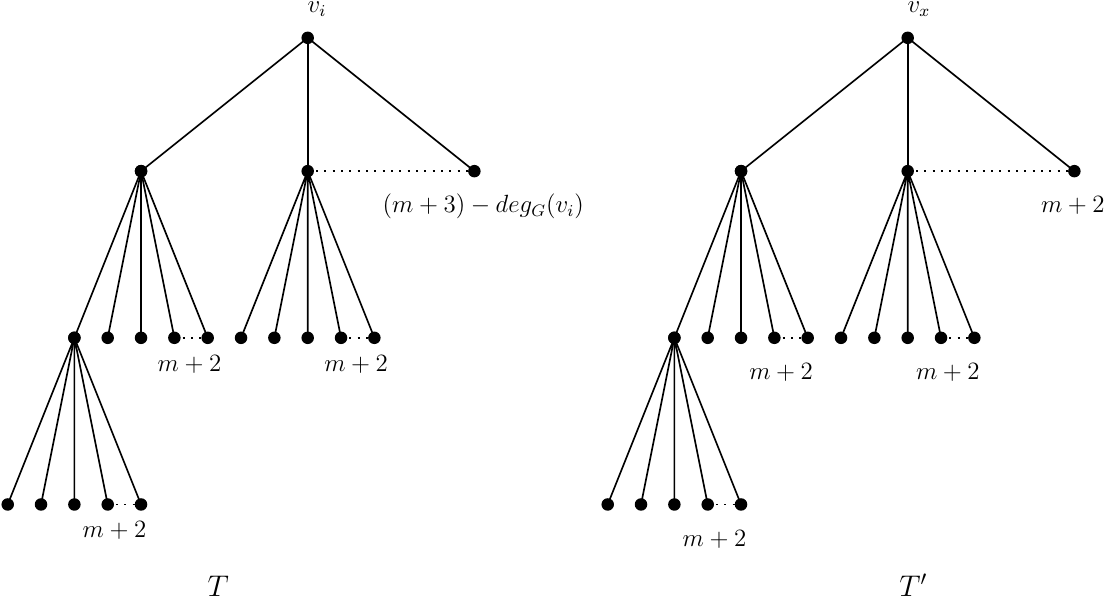}
	\caption{The trees $T$ and $T'$}
	\label{fig:st_np_chordal_tree}
\end{figure}

Given an instance of P$3$CDP, say $G=(V, E)$, we construct an instance of MSTDP. The construction is as follows: let $V=\{v_1, v_2, \ldots, v_n\}$ and $E= \{e_1, e_2, \ldots, e_m\}$.
\begin{itemize}
	\item[$1$]  For each vertex $v_i\in V$, we consider a tree $T$ (shown in Figure \ref{fig:st_np_chordal_tree}) with $v_i$ as the root and degree of the root is $(m+3)-deg_{G}(v_i)$. Also, for each edge, $e_j\in E$, we consider a vertex $v_{e_i}$ and consider another tree $T'$ (shown in Figure \ref{fig:st_np_chordal_tree}) with $v_{e_i}$ as the root, where the degree of the root is $m+2$.
	
	\item[$2$] For each edge $e_j\in E$, we take another vertex $e_j$ in $G'$ and also take another extra vertex $e$ in $G'$. Let $A=\{e_1,e_2,\ldots ,e_m,e\}$. We make a complete graph with vertex set $A$.
	
	\item[$3$] We take another extra three vertices $v_a$, $v_e$ and $v_b$ and consider three trees $T'$ (shown in Figure \ref{fig:st_np_chordal_tree}) with $v_a$, $v_e$ and $v_b$ as the roots, respectively.
	
	\item[$5$] Next we add the following edges: for every edge $e_k=(v_i,v_j)\in E$, we join the edges $(e_k, v_i)$, $(e_k, v_j)$, $(e_k, v_{e_k})$. Also we add the edges $(e, v_a)$, $(e, v_e)$, $(e, v_b)$.
	
	\item[$6$] Finally, we set $k=m+4$.
	
\end{itemize}

\begin{figure}[htbp!]
	\centering
	\includegraphics[scale=0.60]{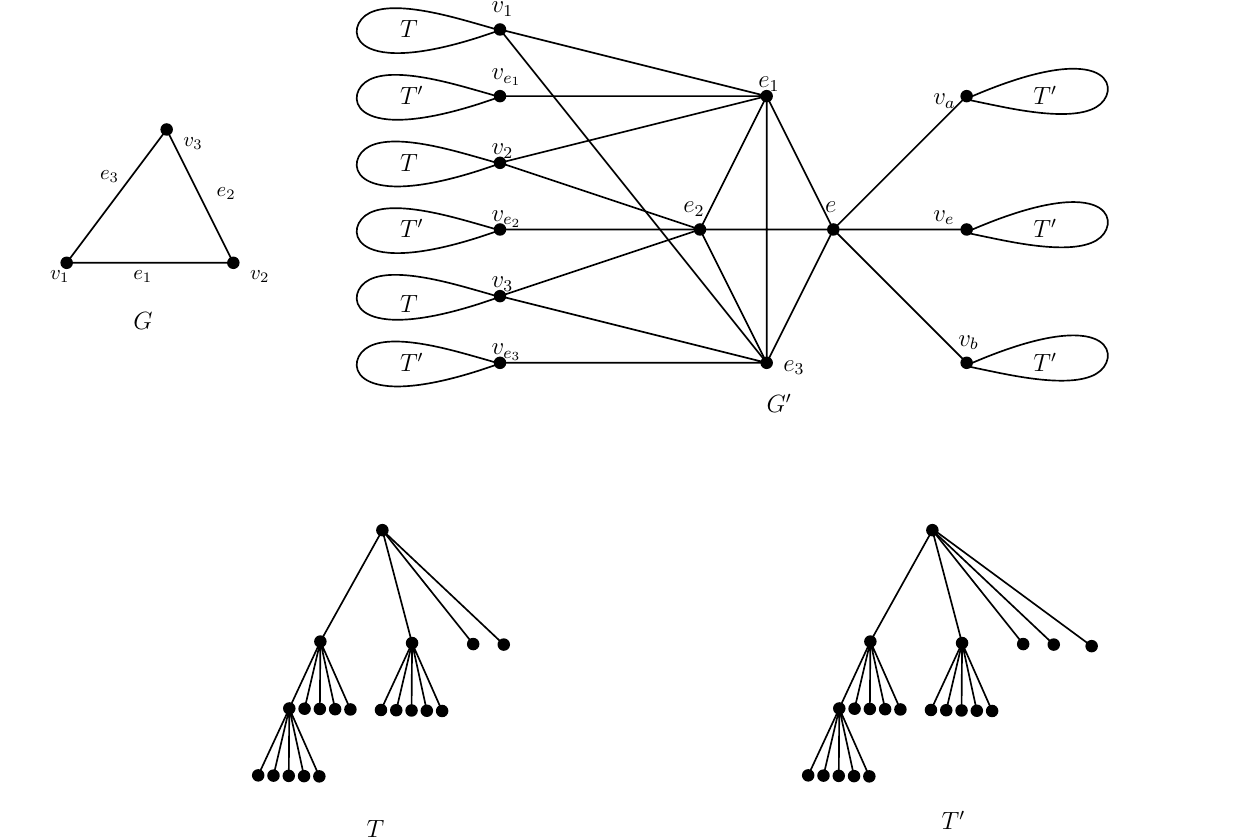}
	\caption{Construction of $G'$ from $G$}
	\label{fig:chordalnp_strong_transitivity}
	
\end{figure}

Note that $G'$ is a chordal graph. The construction from $G$ to $G'$ is illustrated in Figure \ref{fig:chordalnp_strong_transitivity}.

Next, we show that  $G$ has a proper 3-coloring if and only if $G'$ has a strong transitive partition of size $k$. For the forward direction, we have the following lemma.

\begin{lem}
	If $G=(V,E)$ has a proper 3-coloring, then $G'=(V', E')$ has a strong transitive partition of size $k$.
\end{lem}
\begin{proof}
	Given a proper 3-coloring $g$ from $V$ to $\{1,2,3\}$, a strong transitive partition of size $k$, say $\pi=\{V_1,V_2, \ldots,V_k\}$ can be obtain in the following ways:

	\begin{figure}[htbp!]
		\centering
		\includegraphics[scale=0.65]{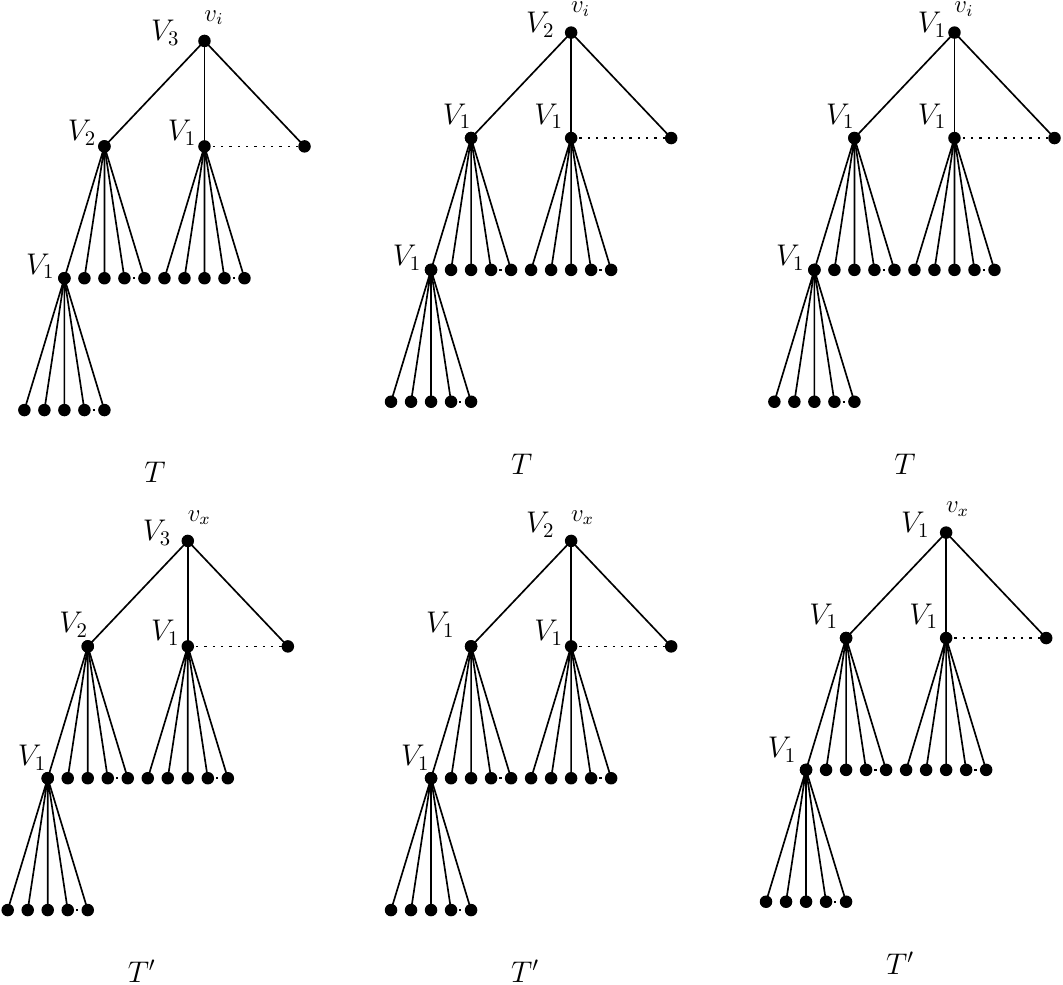}
		\caption{Partition of $T$ and $T'$ into $V_1, V_2$ and $V_3$. All the leaves are in $V_1$}
		\label{fig:st_np_chordal_tree_coloring}
	\end{figure}

	\begin{enumerate}
		\item  If $g(v_i)=q$, then $v_i \in V_q$, for all  $v_i\in V(G)$.
		
		\item  $v_a \in V_3$, $v_{e} \in V_2$ and $v_b \in V_1$. 
		
		\item For each $v_{e_j}$ vertex corresponding an edge $e_j$ with end points $v_x$ and $v_y$ in $G$, assign $v_{e_j} \in V_l$, where $l= \{1, 2, 3\} \setminus \{g(v_x),g(v_y)\}$. Put the other vertices of the trees $T$ and $T'$ in $V_1, V_2$ and $V_3$ based on their root. This is illustrated in Figure \ref{fig:st_np_chordal_tree_coloring}. 
		
		\item Let $e_j\in V_{3+j}$, $1\leq j\leq m+3$, and $e\in V_{m+4}$. 
		
	\end{enumerate}

	Let $H$ be the complete graph induced by $A$. Since $H$ is a complete graph, then $V_i$ strongly dominates $V_j$ for $4\leq i<j\leq k$. Also, for each $i=1, 2, 3$, every vertex of $A$ is adjacent to a vertex of $V_i$, and the degree of that vertex is equal to the degree of a vertex of $A$. Therefore, for each $i=1, 2, 3$, $V_i$ strong dominates $V_j$ for all $j>3$. At the end, from Figure \ref{fig:st_np_chordal_tree_coloring}, it is clear that $V_i$ strongly dominates $V_j$ for $1\leq i<j\leq 3$. Hence, $\pi$ is a strong transitive partition of $G'$ of size $k$. Therefore, if $G$ has a proper 3-coloring, then $G'$ has a strong transitive partition of size $k$.
\end{proof}

Next, we show the converse of the statement. For this, we first prove the following claim. 
\begin{cl}\label{claim:123 and rest_strong}
	Let  $\pi=\{V_1,V_2,\ldots ,V_k\}$ be a strong transitive partition of $G'$ of size $k$ such that $|V_k|=1$. Then the sets $V_4, V_5,\ldots V_k$ contain only vertices from $A$ and the sets $V_1, V_2, V_3$ contain only vertices from $V'\setminus A$.
\end{cl}

\begin{proof}

	We divide the proof into two cases:
	
	\begin{ca}\label{st_case_npchordal}
		$e\in V_{m+4}$
	\end{ca}
	
	Since the degree of each vertex from $\{v_a, v_e, v_b\}$ is $m+3$ and adjacent with three vertices having a degree is more than or equal to $m+3$, then they cannot be in $V_p$, $p\geq 5$. Now, if any vertex from $\{v_a, v_e, v_b\}$ is in $V_4$, then $e$ must be in $V_i$, $1\leq i\leq 3$, a contradiction as $e\in V_k$ and $k\geq 4$. Therefore, the vertices from $\{v_a,v_e,v_b\}$ are belong to $V_p, 1\leq p \leq 3$. The vertex $e\in V_k$, to strongly dominate $e$, each set in $\{V_1, V_2, \ldots, V_{m+3}\}$ must contain at least one vertex from $N_G'(e)=\{e_1, e_2, \ldots, e_m, v_a, v_e, v_b\}$. Since $e$ is adjacent with exactly $m+3$ vertices and the set $N_G'(e)=\{e_1, e_2, \ldots, e_m, v_a, v_e, v_b\}$ contains exactly $m+3$ vertices, each $V_i, 1\leq i \leq m+3$ contains exactly one vertex from $N_G'(e)$. Since $\{v_a, v_e, v_b\}$ belong to $V_p$ for some $1\leq p \leq 3$, it follows that vertices from $\{e_1, e_2, \ldots, e_m\}$ belong to $V_p$, $p\geq 4$. Hence, the vertices of $A$ belong to $\{V_4, V_5,\ldots V_{m+4}\}$. Note that none of the vertices from $\{v_1, v_2, \ldots, v_n, v_{e_1}, v_{e_2}, \ldots, v_{e_m}\}$ belong to $V_p$ for some $p\geq 4$. Because otherwise, there exists a vertex of $A$ must be in $V_3$. But this contradicts the fact that the vertices of $A$ belong to $\{V_4, V_5,\ldots V_k\}$. Since the number of neighbours having more than or equal degree of every other vertices is at most $2$, they cannot belong to $V_p$, $p\geq 4$. Therefore, $V_4, V_5,\ldots V_{m+4}$ contain only vertices from $A$ and $V_1, V_2, V_3$ contain only vertices from $V'\setminus A$.

	\begin{ca}
		$e\notin V_{m+4}$
	\end{ca}
	
	Since $\pi$ is a strong transitive partition, for any $x\in V_{m+4}$, $deg(x)\geq m+3$ and has at least $m+3$ neighbour with a degree at least $deg(x)$. As $deg(v_i)=m+3$ and $v_i$ has at most $m+2$ neighbour having degree at least $m+3$, so $v_i$ cannot be in $V_{m+4}$. Similarly, we can prove that any vertex other than $\{e_1, e_2, \ldots, e_m, e\}$ cannot belong to $V_{m+4}$. Since, $e\notin V_{m+4}$, without loss of generality assume $e_1\in V_{m+4}$, where $e_1$ is the vertex of $G'$ corresponding to the edge $e_1=v_1v_2\in E$. Now we show that $v_1$ and $v_2$ belong to the first three sets in $\pi$. Let $v_1\in V_l$ and $v_2\in V_t$, where $t\leq l$. If possible, let $l\geq 4$. Since $e_1\in V_{m+4}$, to dominate $e_1$, each set in $\{V_1, V_2, \ldots, V_{m+4}\}$ must contain at least one vertex from $N_{G'}(e_1)=\{e_2, e_3, \ldots, e_m, e, v_1, v_2, v_{e_1}\}$. Since $e_1$ is adjacent with exactly $m+3$ vertices, each $V_i$, $1\leq i \leq m+3$ contains exactly one vertex from $N_{G'}(e_1)$. So, if $l\geq 4$, then to strongly dominate $v_1$, each set in $\{V_3, V_4, \ldots, V_{l-1}\}$ contains exactly one vertex  from $\{e_2, e_3, \ldots, e_m\}$. Now, to strongly dominate $e_1$, each set in  $\{V_{l+1}, V_{l+2}, \ldots, V_{m+3}\}$ contains exactly one vertex from $\{e_2, e_3, \ldots, e_m, e\}$. The vertex $e$ cannot belong to $V_q, q\geq l+1$, because $V_l$ contains $v_1$ and any $\{v_a, v_e, v_b\}$  cannot be in $V_l$ as $l\geq 4$.  Also $e$ cannot belong to $V_p, 3\leq p \leq l-1$, as $v_1\in V_l$ and each of $\{V_3, V_4, \ldots, V_{l-1}\}$ contains exactly one vertex from $\{e_2, \ldots, e_m\}$. Therefore, the vertex $e$  must be in  $V_i, 1\leq i\leq 2$. So each of $\{V_3, V_4, \ldots, V_{l-1},V_{l+1}, \ldots, V_{m+3}\}$ contains exactly one vertex from  $\{e_2, \ldots, e_m\}$ only. But the number of vertices in $\{e_2, \ldots, e_m\}$ is $m-1$ whereas we need $k-4=m$ vertices. Therefore, $l$ cannot be more than $3$. Note that $v_{e_1}$ cannot be in $V_j$ for $j\geq 4$ as  no of neighbour of $v_{e_1}$ having degree at least $deg(V_{e_1})=m+3$ other than $e_1$ is $2$. Therefore, the vertices $\{v_1, v_{e_1}, v_2\}$ belong to $V_p$ for $1\leq p \leq 3$. With similar arguments as in Case \ref{st_case_npchordal}, we can say that the vertices of $A$ belong to $\{V_4, V_5,\ldots V_k\}$. We can further claim that, the vertices of $\{v_1, v_2, \ldots, v_n, v_{e_1}, v_{e_2}, \ldots, v_{e_m}\}$ belong to $V_p$ for $1\leq p \leq  3$ and the other vertices of $G'$ belong to $V_p$ for $1\leq p \leq  3$. Therefore, $V_4, V_5,\ldots V_k$ contain only vertices from $A$ and the sets $V_1, V_2, V_3$ contain only vertices from $V'\setminus A$.
\end{proof}

Using the claim, we show that $G$ has a proper $3$-coloring.

\begin{lem}
	If $G'$ has a strong transitive partition of size $k$, then $G$ has a proper $3$-coloring.
\end{lem}
\begin{proof}
	Let $\pi=\{V_1,V_2,\ldots ,V_k\}$ be a strong transitive partition of $G'$ of size $k$. Since $\pi$ is also a transitive partition, from \cite{haynes2019transitivity} we can assume that $|V_k|=1$. Let us define a coloring of $G$, say $g$, by labelling $v_i$ with color $p$ if its corresponding vertex $v_i$ is in $V_p$. The previous claim ensures that $g$ is a $3$-coloring. Now we show that $g$ is a proper coloring. Let $e_t=v_iv_j\in E$ and let its corresponding vertex $e_t$ in $G'$ belong to some set $V_p$ with $p\geq 4$. This implies that the vertices $\{v_i, v_j, v_{e_t}\}$ must belong to different sets from $V_1, V_2, V_3$. Therefore, $g(v_i)\neq g(v_j)$ and hence $g$ is a proper coloring of $G$.
\end{proof}

Therefore, we have the following main theorem of this section:

\begin{theo}
	The \textsc{MSTDP} is NP-complete for chordal graphs.
\end{theo}

\section{Linear-time algorithms}

\subsection{Trees}

In this subsection, we design a linear-time algorithm for finding the strong transitivity of a given tree $T=(V, E)$. We design our algorithm in a similar way to the algorithm for finding the Grundy number of an input tree presented in \cite{hedetniemi1982linear}. First, we give a comprehensive description of our proposed algorithm.

\subsubsection{Description of the algorithm:}
Let $T^c$ denote a rooted tree rooted at a vertex $c$ and $T_v^c$ denote the subtree of $T^c$ rooted at a vertex $v$. With a small abuse of notation, we use $T^c$ to denote both the rooted tree and the underlying tree. To find the strong transitivity of $T=(V, E)$, we first define the \emph{strong transitive number} of a vertex $v$ in $T$. The strong transitive number of a vertex $v$ in $T$ is the maximum integer $p$ such that $v\in V_p$ in a strong transitive partition $\pi=\{V_1, V_2, \ldots, V_k\}$, where the maximum is taken over all strong transitive partition of $T$. We denote the strong transitive number of a vertex $v$ in $T$ by $st(v, T)$. Note that the strong transitivity of $T$ is the maximum strong transitive number that a vertex can have; that is, $Tr_{st}(T)=\max\limits_{v\in V}\{st(v, T)\}$. Therefore, our goal is to find a strong transitive number of every vertex in the tree. Now we define another parameter, namely the \emph{rooted strong transitive number}. The \emph{rooted strong transitive number} of $v$ in $T^c$ is the strong transitive number of $v$ in the tree $T_v^c$ and it is denoted by $st^r(v, T^c)$. Therefore,  $st^r(v, T^c)= st(v, T_v^c)$. To this end, we define another parameter, namely the \emph{modified rooted strong transitive number}. The \emph{modified rooted strong transitive number} of $v$ in $T^c$ is the strong transitive number of $v$ in the tree $T_v^c$, considering the $deg(v)$ as $deg(v)+1$ if $v$ is a non-root vertex and for root vertex $deg(v)$ as $deg(v)$. We denote it by $mst^r(v, T^c)$. Note that the value modified rooted strong transitive number of a vertex depends on the rooted tree, whereas the strong transitive number is independent of the rooted tree. Also, for the root vertex $c$, $st^r(c, T^c)= mst^r(c, T^c)=st(c,T)$. We recursively compute the modified rooted strong transitive number of the vertices of $T^c$ in a bottom-up approach. First, we consider a vertex ordering $\sigma$, which is the reverse of BFS ordering of $T^c$. For a leaf vertex $c_i$, we set $mst^r(c_i, T^c)=1$. For a non-leaf vertex $c_i$, we call the function \textsc{Strong\_Transitive\_Number$()$}, which takes the modified rooted strong transitive number of children of $c_i$ in $T^c$ as input and returns the modified rooted strong transitive number of $c_i$ in $T^c$. At the end of the bottom-up approach, we have the modified rooted strong transitive number of $c_i$ in $T^c$, that is, $mst^r(c, T^c)$, which is the same as the strong transitive number of $c$ in $T$, that is, $st(c, T)$. After the bottom-up approach, we have the strong transitive number of the root vertex $c$, and the modified rooted strong transitive number of every other vertices in $T^c$. 

Next, we compute the strong transitive number of every other vertex. For a vertex $c_i$, other than $c$, we compute the strong transitive number using the function \textsc{Strong\_Transitive\_Number$()$}, which takes the modified rooted strong transitive number of children of $c_i$ in $T^{c_i}$ as input. Let $y$ be the parent of $c_i$ in $T^c$. Note that, except $y$, the modified rooted strong transitive number of children of $c_i$ in $T^{c_i}$ is the same as the modified rooted strong transitive number in $T^c$. We only need to compute the modified rooted strong transitive number of $y$ in $T^{c_i}$. We use another function called \textsc{Strong\_Mark\_Required$()$} for this. This function takes a strong transitive number of a vertex $x$ and modified rooted strong transitive number of its children in $T^x$ as input and marks the status of whether a child, say $v$, is required or not to achieve the strong transitive number of $x$. We mark $R(v)=1$ if the child $v$ is required, otherwise $R(v)=0$. We compute the strong transitive number of every vertex, other than $c$, by processing the vertices in the reverse order of $\sigma$, that is, in a top-down approach in $T^c$. While processing the vertex $c_i$, first based on the status marked by \textsc{Strong\_Mark\_Required$()$} function, we calculate the modified rooted strong transitive number of $p(c_i)$ in $T^{c_i}$, where $p(c_i)$ is the parent of $c_i$ in the rooted tree $T^c$. Then, we call \textsc{Strong\_Transitive\_Number$()$} to calculate the strong transitive number of $c_i$. Next, we call the \textsc{Strong\_Mark\_Required$()$} to mark the status of the children, which will be used in subsequent iterations. At the end of this top-down approach, we have a strong transitive number of all the vertices and hence the strong transitivity of the tree $T$. The process of finding $Tr_{st}(T)$ is described in Algorithm \ref{Algo:strong_trasitivity(T)}.

\begin{algorithm}[h]
	
	\caption{\textsc{Strong\_Transitivity(T)} } \label{Algo:strong_trasitivity(T)}
	
	\textbf{Input:} A tree $T=(V, E)$.
	
	\textbf{Output:} Strong transitivity of $T$.
	
	\begin{algorithmic}[1]
		
		%	\State  \textbf{Input:} A tree $T=(V, E)$.
		
		%	\State  \textbf{Output:} Strong transitivity of $T$.
		
		\State  Let $\sigma=(c_1, c_2, \ldots, c_k=c)$ be the reverse BFS ordering of the  vertices of $T^c$, rooted at a vertex $c$. 
		
		\ForAll {$c_i$ in $\sigma$}
		
		\If {$c_i$ is a leaf }
		
		\State $mst^r(c_i, T^c)=1$.

		\Else
		
		\State $mst^r(c_i, T^c)$ = \textsc{Strong\_Transitive\_Number}$(mst^r(c_{i_1},T^c), \ldots , mst^r(c_{i_k}, T^c) )$.
		
		~~~~~~~ /*where $c_{i_1}, c_{i_2}, \ldots, c_{i_k}$ be the children of $c_i$ and $mst^r(c_{i_1}, T^c)\leq mst^r(c_{i_2}, T^c)\leq \ldots \leq mst^r(c_{i_k}, T^c)$*/
		
		\EndIf
		
		\EndFor

		\State \textsc{Strong\_Mark\_Required}$(st(c, T), mst^r(c_{1}, T^c),  \ldots, mst^r(c_{k}, T^c) )$
		
		~~~~~~~ /*where $c_{1}, c_{2}, \ldots, c_{k}$ be the children of $c$ and $mst^r(c_{1}, T^c)\leq mst^r(c_{2}, T^c)\leq \ldots \leq mst^r(c_{k}, T^c)$*/

		\ForAll{$c_i\in \sigma'\setminus\{c\}$ } %~~~~~~~~~~~~~~~~~~~~~~~~~~~~~~~~~~~~~~~~~~~~~~~~~~~~~/*where $\sigma'$ is the BFS ordering of $T^c$*/
		
		\If {$R(c_i)=0$}
		
		\State $mst^r(p(c_i), T^{c_i}) = st(c_i, T)$
		
		\Else
		
		\State $mst^r(p(c_i), T^{c_i}) = st(c_i, T) -1$
		
		\EndIf

		\State $st(c_i, T) =$\textsc{strong\_Transitive\_Number}$(mst^r(c_{i_1}, T^{c_i}),  \ldots , mst^r(c_{i_k}, T^{c_i}) )$.

		\State \textsc{Strong\_Mark\_Required}$(st(c_i, T), mst^r(c_{i_1}, T^{c_i}), \ldots , mst^r(c_{i_k}, T^{c_i}) )$

		~~/*where $c_{i_1}, c_{i_2}, \ldots, c_{i_k}$ be the neighbours of $c_i$ and $mst^r(c_{i_1}, T^c)\leq mst^r(c_{i_2}, T^c)\leq \ldots \leq mst^r(c_{i_k}, T^c)$*/
		
		~~~~~~~~~~~~~~~~~~~~~~~~~~~~~~~~~~~~~~~~/* Note that $mst^r(c_{i_j}, T^{c_i})=mst^r(c_{i_j}, T^{c})$ if $c_{i_j}$ is a child of $c_i$ in $T^c$*/
		
		\EndFor

		\State $Tr_{st}(T)=\max\limits_{x\in V}\{st(x, T)\}$.

	\end{algorithmic}
	
\end{algorithm}

\subsubsection{Proof of correctness}

In this subsection, we give the proof of the correctness of Algorithm \ref{Algo:strong_trasitivity(T)}. It is clear that the correctness of Algorithm \ref{Algo:strong_trasitivity(T)} depends on the correctness of the functions used in the algorithm. First, we show the following two lemmas, which prove the correctness of Strong\_Transitive\_Number$()$ function.

\begin{lem}\label{tree_lemma_strong_transitivity}
	Let $x$ be a child of $T^c$  and $y$ be its parent in $T^c$. Also, let  $mst^r(v,T)=t$. Then there exists a strong transitive partition of $T_x^c$, say $\{V_1, V_2, \ldots, V_i\}$ such that $x\in V_i$, for all $1\leq i\leq t$.
\end{lem}

\begin{proof}
	Since $mst^r(x,T^c)=t$, there exists a strong transitive partition $\pi=\{U_1,U_2,\ldots, U_t\}$ of $T_x^c$ such that $x\in U_t$. 
	For each $1\leq i\leq t$, let us define another strong transitive partition $\pi'=\{V_1,V_2,\ldots,V_i\}$ of $T_x^c$ as follows: $V_j=U_j$ for all $1\leq j\leq (i-1)$ and $V_i= \displaystyle{\bigcup_{j=i}^{t}U_j}$. Clearly, $\pi'$ is a strong transitive partition of $T_x^c$ of size $i$ such that $x\in V_{i}$. Hence, the lemma follows.
\end{proof}

\begin{lem}\label{tree_theorem_strong_transitivity}
	Let $v_1, v_2, \ldots, v_k$ are the children of $x$ in a rooted tree $T^c$ and $y$ be the parent of $x$ in  $T^c$. Also, let for each $1\leq i\leq k$, $l_i$ denote the modified rooted strong transitive number of $v_i$ in $T^c$ with $l_1\leq l_2\leq \ldots\leq l_k$ and $p(y)=1$ when $y$ exists in $T^c$ otherwise $p(y)=0$. Let $z$ be the largest integer such that there exists a subsequence of $\{l_i: 1\leq i\leq k\}$, say $(l_{i_1}\leq l_{i_2}\leq \ldots \leq l_{i_{z}})$ such that $l_{i_{p}}\geq p$, for all $1\leq p\leq z$ and $deg(v_{i_j})\geq deg(x)+p(y)$. Then, the modified rooted strong transitive number of $x$ in the underlying tree $T^c$ is $1+z$, that is, $mst^r(x, T^c)=1+z$.
\end{lem}

\begin{proof}
	For each $1\leq j\leq z$, let us consider the subtrees $T_{v_{i_j}}^c$. It is also given that $mst^r(v_{i_j},T^c)=l_{i_j}$, for $j\in \{1, 2, \ldots, {z}\}$. For all $1\leq p\leq z$, since $l_{i_{p}}\geq p$,  by  Lemma \ref{tree_lemma_strong_transitivity}, we know that there exists strong transitive partitions $\pi^{p}=\{V_1^{p}, V_2^{p}, \ldots, V_p^{p}\}$ of  $T_{v_{i_{p}}}^c$ such that $v_{i_{p}}\in V_p^{p}$ and $deg(v_{i_p})\geq deg(x)+1$. Let us consider the partition of $\pi=\{V_1, V_2, \ldots, V_z, V_{z+1}\}$ of $T_x^c$ as follows: $V_i=\displaystyle{\bigcup_{j=i}^{z}V_i^j}$, for $\leq i\leq z$, $V_{z+1}=\{x\}$ and every other vertices of $T$ are put in $V_1$. Clearly, $\pi$ is a strong transitive partition of $T_x^c$. Also it is given that $deg(v_{i_j})\geq deg(x)+1$. Therefore, $mst^r(x,T^c)\geq 1+z$.
	
	Next, we show that $mst(x, T^c)$ cannot be more than $1+z$. If possible, let $mst(x,T^c)\geq 2+z$. Then by Lemma \ref{tree_lemma_strong_transitivity}, we have that there exists a strong transitive partitions $\pi=\{V_1, V_2, \ldots, V_{2+z}\}$ such that $x\in V_{2+z}$. This implies that for each $1\leq i\leq 1+z$, $V_i$ contains a neighbour of $x$, say $v_i$ such that the modified rooted strong transitive number of both $v_i$ is greater or equal to $i$, that is, $l_i\geq i$ and $deg(v_i)\geq deg(x)+1$. The set $\{l_i | 1\leq i\leq 1+z\}$ forms a desired subsequence of $\{l_i: 1\leq i\leq k\}$, contradicting the maximality of $z$. Hence, $mst(x,T)=1+z$.
\end{proof}

Note that in line $6$ of Algorithm \ref{Algo:strong_trasitivity(T)}, when Strong\_Transitive\_Number$()$ is called, then it returns the strong transitive number of $c_i$ in $T^c_{c_i}$ which is in fact the modified rooted strong transitive number of $c_i$ in $T^c$. And in line $13$ of Algorithm \ref{Algo:strong_trasitivity(T)}, when Strong\_Transitive\_Number$()$ is called, then it returns the strong transitive number of $c_i$ in $T^{c_i}$ which is same as $st(c_i, T)$. From Lemma \ref{tree_lemma_strong_transitivity} and \ref{tree_theorem_strong_transitivity}, we have the function \textsc{Strong\_Transitive\_Number$()$}.

\begin{algorithm}[h]
	
	\caption{\textsc{Strong\_Transitive\_Number$(mst^r(v_{1}, T^x), mst^r(v_{2}, T^x), \ldots, mst^r(v_{k}, T^c) )$} } \label{Algo:strongtransNumber(x)}

	\textbf{Input:} Modified rooted strong transitive numbers of children of $x$ in the tree $T^c$ with $mst^r(v_{1}, T^x)\leq mst^r(v_{2}, T^x)\leq \ldots \leq mst^r(v_{k}, T^x)$ and $p(y)$, where $y$ is the parent of $x$ in $T^c$.
	
	\textbf{Output:} Modified rooted strong transitive number of $x$ in the underlying tree $T^c$, that is, $mst^r(x, T^c)$.
	
	\begin{algorithmic}[1]
		
		\State $mst^r(x, T^c)$ $\leftarrow$ 1.
		
		\ForAll {$i\leftarrow 1$ to $k$}
		
		\If {$mst^r(v_i,T^c)\geq mst^r(x, T^c)$ and $deg(v_i)\geq deg(x)+p(y)$}
		
		\State $mst^r(x, T)=mst^r(x, T)+1$.
		
		\Else
		
		\State  $mst^r(x, T^c)=mst^r(x, T^c)$
		
		\EndIf
		
		\EndFor

		\State \Return($mst^r(x, T^c)$).
		
	\end{algorithmic}
	
\end{algorithm}

Next, we prove the correctness of Strong\_Mark\_Required$()$. Let $T^x$ be a rooted tree and $st(x, T)=z$. A child $v$ of $x$ is said to be required if the $st(x, T^x \setminus T_v^x)=z-1$. The function returns the required status of every child of $x$ by marking $R(v)=1$ if it is required and $R(v)=0$ otherwise. The children of $x$ that are required can be identified using the following lemma.

\begin{lem} \label{strong_transitivity_required_descendants}
	Let $T^x$ be a tree rooted at $x$ and $v_1, v_2, \ldots, v_k$ be its children in $T^x$. Also, let the strong transitive number of $x$ be $z$ and for each $1\leq i\leq k$, let $l_i$ denote the modified rooted strong transitive number of $v_i$ in $T^x$. Moreover, let  for all $1\leq i\leq p$, $deg(v_i)<k$ and for all $p+1\leq i\leq k$, $deg(v_i)\geq k$ and $l_{p+1}\leq l_{p+2}\leq \ldots\leq l_k$. Then the following hold:
	
	\begin{enumerate}
		\item[(a)] If $k=z-1$, then $R(v_i)=1$ for all $1\leq i\leq k$.
		
		\item[(b)] Let $k>z-1$. If $k-p=z-1$, then for all $1\leq i\leq p$, $R(v_i)=0$ and for all $p+1\leq i\leq k, R(v_i)=1$.
		
		\item[(c)] Let $k>z-1$ and $k-p>z-1$. Then for all $1\leq i\leq p$, $R(v_i)=0$ and for all $p+1\leq i\leq k-z+1, R(v_i)=0$.
		
		\item[(d)] Let $k>z-1$ and $k-p>z-1$. Also, let $k-z+2\leq i\leq k$. If for all $j$, $k-z+2\leq j\leq i$, $l_{j-1}\geq j-(k-z+1)$  then $R(v_i)=0$. 
		
		\item[(e)]Let $k>z-1$, $k-p>z-1$ and $k-z+2\leq i\leq k$. If there exists $j$ in $k-z+2\leq j\leq i$ such that $l_{j-1}< j-(k-z+1)$ or then $R(v_i)=1$.
		
	\end{enumerate}
\end{lem}

\begin{proof}
	
	Note that $k-p\geq z-1$ as $st(x, T^x)=z$ and $deg(x)=k$.
	
	$(a)$ Since $st(x, T^x)=z$, by the Lemma \ref{tree_lemma_strong_transitivity}, there exists a strong transitive partition of $T$, say $\pi=\{V_1, V_2, \ldots, V_{z}\}$, such that $x\in V_z$. In that case, all the vertices in $\{v_1, v_2, \ldots, v_k\}$ must be in $V_1, V_2, \ldots, V_{z-1}$ and each set $V_i$ contains at least one of these vertices. Since $k=z-1$, each set $V_i$ contains exactly one vertex from $\{v_1, v_2, \ldots, v_k\}$. Therefore, if we remove any $v_i$ from the tree, the strong transitive number of $x$ will decrease by $1$. Hence, every $v_i$ is required, that is, $R(v_i)=1$ for all $1\leq i\leq k$.
	
	$(b)$ Since $st(x, T^x)=z$, by the Lemma \ref{tree_lemma_strong_transitivity}, there exists a strong transitive partition of $T^x$, say $\pi=\{V_1, V_2, \ldots, V_{z}\}$, such that $x\in V_z$. Since the vertices from $\{v_1, v_2, \ldots, v_p\}$ have degree less than $x$, they are not use to strong dominate $x$. Therefore,  if we remove any $v_i$ ($1\leq i\leq p$) from the tree, then the strong transitive number of $x$ will be unchanged. Hence, for each $1\leq i\leq p$,  $v_i$ is not required, that is, $R(v_i)=0$ for all $1\leq i\leq p$. So, all the vertices in $\{v_{p+1}, v_{p+2}, \ldots, v_k\}$ must be in $V_1, V_2, \ldots, V_{z-1}$ and each set $V_i$ contains at least one of these vertices. Since $k-p=z-1$, each set $V_i$ contains exactly one vertex from $\{v_{p+1}, v_{p+2}, \ldots, v_k\}$. Therefore, removing any $v_i$ ($p+1\leq i\leq k$) from the tree will decrease the strong transitive number of $x$ by $1$. Hence, every $v_i$ is required, that is, $R(v_i)=1$ for all $p+1\leq i\leq k$.
	
	$(c)$ Let $\pi=\{V_1, V_2, \ldots, V_{z}\}$ be a strong transitive partition of $T^x$ such that $x\in V_z$. As before, the vertices from $\{v_1, v_2, \ldots, v_p\}$ are not required. Hence, $R(v_i)=0$ for all $1\leq i\leq p$.
	Now, in this partition, at least one vertex from $\{v_{p+1}, v_{p+2}, \ldots, v_k\}$  must be in each $V_i$ for $1\leq i\leq z-1$. As the vertices are arranged in increasing order of their modified rooted strong transitive number, without loss of generality, we can assume that $\{v_{p+1}, v_{p+2}, \ldots, v_{k-z+1}\}\subset V_1$ and $v_i\in V_{I_i}$ for each $k-z+2\leq i \leq k$, where $I_i=i-(k-z+1)$. Clearly, if we remove any $v_i$ ($p+1\leq i\leq k-z+1$) from the tree, then the strong transitive number of $x$ will be unchanged. Hence, for each $p+1\leq i\leq k-z+1$, $v_i$ is not required, that is, $R(v_i)=0$ for all $p+1\leq i\leq k-z+1$.

	$(d)$ Let us consider the same strong transitive partition $\pi$ of $T^x$ as in case $(c)$. Let for some $k-z+2\leq i \leq k$,  $v_i$ be a vertex such that $l_{j-1}\geq I_j$ for all $k-z+2\leq j\leq i$, where $I_j= j-(k-z+1)$. We can modify $\pi$ to get a strong transitive partition of $T^x\setminus \{v_i\}$ of size $z$. The modification is as follows: for each $j\in \{k-z+2, k-z+3, \ldots, i\}$, we put $v_{i-1}\in V_{I_i}$ and remove the vertices that are not in $T^x\setminus \{v_i\}$. Therefore, $v_i$ is not required and $R(v_i)=0$ for such vertices.

	$(e)$ Let for some $i$, $k-z+2\leq i \leq k$,  $v_i$ be a vertex such that $l_{j-1}< I_j$ for some $k-z+2\leq j\leq i$, where $I_j= j-(k-z+1)$. Since the modified rooted strong transitive numbers are arranged in increasing order, $l_q< I_j$ for all $p+1\leq q\leq j-1$. Suppose, after deleting the vertex $v_i$, let the strong transitive number of $x$ in $T^x\setminus \{v_i\}$ remain $z$. Let $\pi'=\{V_1, V_2, \ldots, V_{z}\}$ be a strong transitive partition of $T^x\setminus \{v_i\}$ of size $z$ such that $x\in V_{z}$. Since $l_q< I_j$ for all $p+1\leq q\leq j-1$, none of the vertices of $\{v_{p+1}, v_{p+2}, \ldots, v_{j-1}\}$ can be in the sets $V_{I_j}, V_{I_{j+1}}, \ldots, V_{I_k}$. On the other hand, the sets $V_{I_j}, V_{I_{j+1}}, \ldots, V_{I_k}$ must contain at least $(k-j+1)$ vertices from $\{v_{p+1}, v_{p+2}, \ldots, v_{i-1}, v_{i+1}, \ldots, v_{k}\}$, as $\pi'$ is a strong transitive partition of $T^x\setminus \{v_i\}$. Therefore, $V_{I_j}, V_{I_{j+1}}, \ldots, V_{I_k}$ contains at least $(k-j+1)$ vertices from $\{v_{j}, v_{j+1}, \ldots, v_{i-1}, v_{i+1}, \ldots, v_{k}\}$. But there are only $(k-j)$ many vertices available. Hence, the strong transitive number of $x$ in $T^x\setminus \{v_i\}$ cannot be $z$. Therefore, $v_i$ is required and $R(v_i)=1$ for such vertices.
\end{proof}
Note that the condition in case $(e)$ is such that if $R(v_i)=1$ for some $i$, then $R(v_j)=1$ for all $i+1\leq j\leq k$. Based on the lemma, we have the function \textsc{Strong\_Mark\_Required$()$}.

\begin{algorithm}[h]
	
	\caption{\textsc{Strong\_Mark\_Required$(st(x, T), mst^r(v_1, T^x), \ldots,  mst^r(v_k,  T^x) )$}}  \label{Algo:strong_required_descendants(x)}
	
	\textbf{Input:} A rooted tree $T^x$, rooted at a vertex $x$ and modified rooted strong transitive number of the children of $x$, such that $deg(v_i)<k$ for $1\leq i\leq p$ and $deg(v_i)\geq k$, for all $i\geq p+1$ also $mst^r(v_{p+1}, T^c)\leq  \ldots \leq mst^r(v_k, T^c)$.
	
	\textbf{Output:} $R(v)$ value of $v$, $v$ is a child of $x$ in $T^x$.
	\begin{algorithmic}[1]

		\If {the number of children is $z-1$, that is $k=z-1$}
		
		\State $R(v)=1$, for all children $v$ of $x$ in $T^x$.
		
		\EndIf
		
		\If {$k> z-1$}
		
		\State For all $1\leq i\leq k-z+1$, $R(v_i)=0$
		
		\EndIf
		
		\ForAll {$i\leftarrow k-z+2$ to $k$}
		
		\If {$l_{i-1}\geq i-(k-z+1) $}
		
		\State $R(v_i)=0$
		
		\Else 
		
		\State  $R(v_i)=1$
		
		\State {\bf break}
		
		\EndIf
		
		\EndFor
		
		\State For all $j> i$, $R(v_j)=1$

	\end{algorithmic}
	
\end{algorithm}

\subsubsection{Complexity Analysis:}

In the function \textsc{Strong\_Transitive\_Number()}, we find the strong transitive number of vertex $x$ based on the modified rooted strong transitive number of its children. We have assumed that the children are sorted according to their modified rooted strong transitive number. Since the for loop in line $2-6$ of \textsc{Strong\_Transitive\_Number()} runs for every child of $x$, this function takes $O(deg(x))$ time. Similarly, \textsc{Strong\_Mark\_Required()} takes a strong transitive number of a vertex $x$ and modified rooted strong transitive number of its children in $T^x$ as input and marks the status of whether a child, say $v$, is required or not to achieve the strong transitive number of $x$. Here also, we have assumed that the children are sorted according to their modified rooted strong transitive number. Clearly, line $1-2$ of \textsc{Strong\_Mark\_Required()} takes $O(deg(x))$. In line $3-4$, we mark the status for a few children without any checking and for each of the remaining vertices, we mark the required status by checking condition in $O(1)$ time. Therefore, \textsc{Strong\_Mark\_Required()} also takes $O(deg(x))$ time. In the main algorithm \textsc{Strong\_Transitivity(T)}, the vertex order mentioned in line $1$ can be found in linear time. Then, in a bottom-up approach, we calculate every vertex's modified rooted strong transitive numbers. For that we are spending $O(deg(c_i))$ for every $c_i\in \sigma$. Note that we must pass the children of $c_i$ in a sorted order to  \textsc{Strong\_Transitive\_Number()}. But as discussed in \cite{hedetniemi1982linear} (algorithm for finding Grundy number of a tree), we do not need to sort all the children based on their modified rooted strong transitive numbers; sorting the children whose modified rooted strong transitive number is less than $deg(c_i)$, is sufficient. We can argue that this can be done in $O(deg(c_i))$ as shown in \cite{hedetniemi1982linear}. Hence, the loop in line $1-6$ takes linear time. Similarly, we conclude that line $8-14$ takes linear time. Therefore, we have the following theorem:

\begin{theo}
	The \textsc{MSTP} can be solved in linear time for trees.
\end{theo}

\subsection{Transitivity in split graphs} \label{STSG}

A graph $G=(V, E)$ is said to be a \emph{split graph} if $V$ can be partitioned into an independent set $S$ and a clique $K$. In this subsection, we prove that the transitivity of a split graph $G$ is $\omega(G)$, where $\omega(G)$ is the size of a maximum clique in $G$. First, we prove that $Tr_{st}(G)\geq \omega(G)$.

\begin{lem}\label{STSGLM1}
	Let $G=(S\cup K, E)$ be a split graph, where $S$ and $K$ are an independent set and a clique of $G$, respectively. Also, assume that $K$ is the maximum clique of $G$, that is, $\omega(G)=|K|$. Then $Tr_{st}(G)\geq \omega(G)$.
\end{lem}
\begin{proof}
	Let $\omega(G)=t$ and $\{v_1, v_2, \ldots, v_t\}$ be the vertices of  a maximum clique. Also, assume $deg(v_1)\geq deg(v_2)\geq \ldots \geq deg(v_t)$. Consider a vertex partition $\pi=\{V_1, V_2, \ldots, V_{t}\}$ of size $\omega(G)$ by considering each $V_i=\{v_i\}$, for $i\geq 2$ and $V_1=V\setminus \{v_2, v_3, \ldots, v_t\}$. Since the vertices $\{v_1, v_2, \ldots, v_t\}$ form a clique and $deg(v_1)\geq deg(v_2)\geq \ldots \geq deg(v_t)$, $V_i$ strongly dominates $V_j$ for all $1\leq i<j\leq t$. Therefore, $\pi$ forms a strong transitive partition of $G$ with size $t$. Hence, $Tr_{st}(G)\geq t=\omega(G)$.
\end{proof}

Next, in the following lemma, we show that $Tr_{st}(G)=\omega(G)$.

\begin{lem}\label{STSGLM2}
	Let $G=(S\cup K, E)$ be a split graph, where $S$ and $K$ are an independent set and a clique of $G$, respectively. Also, assume that $K$ is the maximum clique of $G$, that is, $\omega(G)=|K|$. Then $Tr_{st}(G)= \omega(G)$. 
\end{lem}
\begin{proof}
	From \cite{santra2023transitivity}, we know that $Tr(G)=\omega(G)+1$ if and only if every vertex of $K$ has a neighbour in $S$. Also, we know that $Tr(G)\geq Tr_{st}(G)$. Now we divide our proof into the following two cases:
	
	\begin{ca}
		A vertex $x\in K$ exists, such that $x$ has no neighbour in $S$.
	\end{ca}
	In this case, form \cite{santra2023transitivity}, we know that $Tr(G)=\omega(G)$. As $Tr_{st}(G)\leq Tr(G)$, so $Tr_{st}(G)\leq \omega(G)$. Again, by the  Lemma \ref{STSGLM1}, we have $\omega(G)\leq Tr_{st}(G)$. Therefore, $Tr_{st}(G)= \omega(G)$.

	\begin{ca}
		
		Every vertex of $K$ has a neighbour in $S$.
	\end{ca}			
	In this case we have $Tr(G)=\omega(G)+1$ \cite{santra2023transitivity}. So, $Tr_{st}(G)\leq Tr(G)=\omega(G)+1$. Suppose $Tr_{st}(G)=\omega(G)+1$ and $\pi=\{V_1, V_2, \ldots, V_{\omega(G)+1}\}$ be a strong transitive partition of $G$ with size $\omega(G)+1$. Since $|K|=\omega(G)$, a set in $\pi$ contains only vertices from $S$. Also, note that for any $s\in S$, $deg(s)<deg(x)$, for all $x\in K$ as  $K$, is the maximum clique of $G$, and every vertex of $K$ has a neighbour in $S$. Let $s\in S$ and $s\in V_{\omega(G)+1}$. Then $deg(s)$ is at least $\omega(G)$, which is impossible. So, no vertices from $S$ in $V_{\omega(G)+1}$. Let $x\in K$ and $x\in V_{\omega(G)+1}$. Also, assume $V_i\subseteq S$. Since $\pi$ is a strong transitive partition, $V_i$ strongly dominates $V_{\omega(G)+1}$. That implies $deg(s)\geq deg(x)$ for some $s\in S$. We have a contradiction as $deg(s)<deg(x)$ , for all $s\in S$ and $x\in K$. Therefore, $Tr_{st}(G)$ cannot be $\omega(G)+1$. Hence, $Tr_{st}(G)< \omega(G)+1$. Again, by the  Lemma \ref{STSGLM1}, we have $\omega(G)\leq Tr_{st}(G)$. Therefore, $Tr_{st}(G)= \omega(G)$.

\end{proof}

From Lemma \ref{STSGLM2}, it follows that computing the strong transitivity of a split graph is the same as computing the maximum clique. Note that finding a vertex partition of $V$ into $S$ and $K$, where $S$ and $K$ are an independent set and a clique of $G$, respectively, and  $\omega(G)=|K|$ can be computed in linear time \cite{hammer1981splittance}. Hence, we have the following theorem:

\begin{theo}
	The \textsc{MSTP} can be solved in linear time for split graphs.
\end{theo}

\section{Conclusion}
In this paper, we have introduced the notion of strong transitivity in graphs, which is a variation of transitivity. We have shown that the decision version of this problem is NP-complete for chordal graphs. On the positive side, we have proved that this problem can be solved in linear time for trees and split graphs. It would be interesting to investigate the complexity status of this problem in other graph classes. Designing an approximation algorithm for this problem would be another challenging open problem.

\section*{Acknowledgements:} Subhabrata Paul was supported by the SERB MATRICS Research Grant (No. MTR/2019/000528). The work of Kamal Santra is supported by the Department of Science and Technology (DST) (INSPIRE Fellowship, Ref No: DST/INSPIRE/ 03/2016/000291), Govt. of India.

\bibliographystyle{alpha}
\bibliography{Strong_Transitivity_bibliography}

\end{document}